\documentclass[reqno]{amsart}

\usepackage{amsmath,amssymb,amsthm}

\usepackage{tikz}
\usepackage{caption}
\usepackage{graphicx}
\usepackage{graphics}

\newtheorem*{thm}{Theorem}

\newtheorem*{conjecture}{Conjecture}

\begin{document}

\title[]{A Note on Approximate Hadamard Matrices}

\author[]{Stefan Steinerberger}
\address{Department of Mathematics, University of Washington, Seattle, WA 98195, USA} \email{steinerb@uw.edu}

 \thanks{The author is supported by the NSF (DMS-2123224).}

\begin{abstract}  A Hadamard matrix is a scaled orthogonal matrix with $\pm 1$ entries. Such matrices exist in certain dimensions: the Hadamard conjecture is that such a matrix always exists when $n$ is a multiple of 4. A conjecture attributed to Ryser is that no circulant Hadamard matrices exist when $n > 4$. Recently, Dong and Rudelson proved the existence of \textit{approximate} Hadamard matrices in all dimensions: there exist universal $0< c < C < \infty$ so that for all $n \geq 1$, there is a matrix $A \in \left\{-1,1\right\}^{n \times n}$ satisfying, for all $x \in \mathbb{R}^n$,
$$ c \sqrt{n} \|x\|_2 \leq \|Ax\|_2 \leq C \sqrt{n} \|x\|_2.$$
 We observe that, as a consequence of the existence of flat Littlewood polynomials, circulant approximate Hadamard matrices exist for all $n \geq 1$.
\end{abstract}
\maketitle

\vspace{-5pt}

\section{Introduction and Result}

\subsection{Hadamard matrices.}
Hadamard matrices are $n \times n$ matrices all of whose entries are $\pm 1$ which are rescaled orthogonal matrices: the rows are orthogonal and thus, in particular,
$ \|H x\|_2 = \sqrt{n} \cdot \|x\|_2.$
Small examples of such matrices are
$$ \begin{pmatrix} 1 \end{pmatrix}, \quad \begin{pmatrix} 1 & 1 \\ 1 & -1 \end{pmatrix},  \quad \begin{pmatrix} -1 & 1 & 1 & 1 \\ 1 & -1 & 1 & 1 \\ 1 & 1 & -1 & 1 \\ 1 & 1 & 1 & -1 \end{pmatrix}, \quad \dots$$
Hadamard matrices are central objects in a number of different areas, we refer to the books by Agaian \cite{ag} and Horadam \cite{hor}. Sylvester \cite{sylv} was the first prove existence in dimensions $n = 2^k$ by noting that if $H$ is an $n \times n$ Hadamard matrix, then $H \otimes H_2$ is a $2n \times 2n$ Hadamard matrix. It is known that if $H$ is a $n \times n$ Hadamard matrix and $n \geq 4$, then $n$ needs to be a multiple of 4. The famous Hadamard conjecture, sometimes ascribed to Paley \cite{paley}, states that this necessary condition is also sufficient. A conjecture ascribed to Ryser \cite{ryser}, but possibly older  \cite{schm}, is that the $4 \times 4$ Hadamard matrix shown above is the largest circulant Hadamard matrix.

\subsection{Approximate Hadamard matrices.} Motivated by an explicit problem of Riesz bases in random frames, Dong and Rudelson \cite{dong} recently introduced the (intrinsically interesting) concept of \textit{approximate} Hadamard matrix: this is a matrix $A$ whose entries are $\pm 1$ that is `close' to an isometry in the sense that
$$ c\|x\| \leq \|Ax\| \leq C\|x\| \qquad \mbox{and the ratio}~C/c~\mbox{is small.}$$
$c$ is the smallest singular value of the matrix $A$ while $C$ is the largest. The ratio $\kappa(A) = C/c \geq 1$ is also known as the condition number which is 1 if and only if the matrix is a (scaled) orthogonal matrix. One has $\kappa(A) > 1$ when $n$ is not a multiple of 4 (since no Hadamard matrix exists): can the condition number stay bounded? Can it be close to 1? Other notions of `almost' Hadamard matrices exist \cite{ban, jaming, park}.

\begin{thm}[Dong-Rudelson \cite{dong}]
There exist universal $0< c < C < \infty$ such that for all $n \geq 1$, there exists a matrix $A \in \mathbb{R}^{n \times n}$ whose entries are $\pm 1$ such that
$$ c \sqrt{n} \|x\|_2 \leq \|Ax\|_2 \leq C \sqrt{n} \|x\|_2.$$
\end{thm}

The proof is highly nontrivial and uses a number of sophisticated ingredients. One ingredient is a construction due to Matolcsi-Rusza \cite{mat} to build approximate $q \times q$ Hadamard matrices by using quadratic residues in $\mathbb{Z}_q$: the structure of Gauss sums implies flatness of the Fourier spectrum. Another ingredient is Vinogradov's theorem that every sufficiently large odd number is the sum of three prime numbers, a sophisticated gluing procedure is then used to conclude the result.

\subsection{Result}
The main purpose of this short note is to note that the above result is true under the additional condition that the matrix is a circulant matrix with entries $\pm 1$. Recall that a circulant matrix is a matrix of the form
\[
A = \begin{bmatrix}
    a_0 & a_{n-1} & a_{n-2} & \ldots & a_1 \\
    a_1 & a_0 & a_{n-1} & \ldots & a_2 \\
    \vdots & \vdots & \vdots & \ddots & \vdots \\
    a_{n-1} & a_{n-2} & a_{n-3} & \ldots & a_0
\end{bmatrix}
\]
These matrices have a number of desirable properties: in particular, they are diagonalized by the discrete Fourier Transform and inversion is fast. They can also be interpreted as discrete convolution operators. Hadamard matrices are useful and circulant matrices are useful, their (approximate) combination may also be useful.
\begin{thm}
There exist universal $0< c < C < \infty$ such that for all $n \geq 1$, there exists a \emph{circulant} matrix $A \in \mathbb{R}^{n \times n}$ whose entries are $\pm 1$ such that
$$ c \sqrt{n} \|x\|_2 \leq \|Ax\|_2 \leq C \sqrt{n} \|x\|_2.$$
\end{thm}

The argument is short at the cost of invoking a powerful result: we prove that flat Littlewood polynomials can be used to construct circulants that are also approximately Hadamard. Balister-Bollob\'as-Morris-Sahasrabudhe-Tiba \cite{balister}, solving a long-standing conjecture of Littlewood \cite{little}, proved that flat Littlewood polynomials exist and the result follows. As a consequence, for a suitable (absolute) choice of constants $0 < c < C < \infty$, one might be inclined to believe that a great many approximate Hadamard matrices should exist: even the strong requirement of being a circulant is not prohibitive. It could be nice to have explicit constructions of approximate Hadamard matrices with $C/c$ guaranteed to be small.

\section{Proof}
\begin{proof}
The singular values of a circulant $A$ are given by the absolute value of its (possibly complex) eigenvalues (see, for example, \cite{karner}). There are several ways of seeing this: a canonical approach is to use the Fourier matrix
$$ \mathcal{F} = \frac{1}{\sqrt{n}} \left( e^{-2 \pi i m k/n}\right)_{m,k = 0}^{n-1}$$
and note that every circulant matrix can be written as
$ A = \mathcal{F}^{-1} D \mathcal{F}$
for some diagonal matrix $D$. This is simply one way of stating that (discrete) convolution is diagonalized by the (discrete) Fourier transform. The singular values of $A$ are the square root of the eigenvalues of $A^T A$ which can be computed from that representation. The question is therefore whether we can find a circulant matrix with the property that the absolute value of its eigenvalues are all comparable up to a universal multiplicative factor that is independent of the size of the matrix. An (unordered) list of the eigenvalues $\lambda_j$ of a circulant os given by, where $0 \leq j \leq n-1$, 
$$ \lambda_j = a_0 + a_1 w^j + a_2 w^{2j} + \dots + a_{n-1} w ^{(n-1)j} \qquad \mbox{with} \quad w = \exp(2 \pi i /n)$$
being a primitive $n-$th root of unity.
 Introducing the polynomial
$$ p(z) = a_0 + a_1 z + a_2 z^2 + \dots + a_{n-1} z^{n-1}$$
we see that
\begin{align*}
 \kappa(A) = \frac{\max_{\|x\|=1} \|Ax\|}{\min_{\|x\|=1} \|Ax\|} &=\frac{\max_{0 \leq j \leq n} |\lambda_j|}{\min_{0 \leq j \leq n} |\lambda_j|} \\
&= \frac{\max_{z^n = 1} |p(z)|}{\min_{z^n=1}|p(z)|} \leq  \frac{\max_{|z| = 1} |p(z)|}{\min_{|z| = 1}|p(z)|}.
\end{align*}
A polynomial 
$$ p(z) = a_0 + a_1 z + a_2 z^2 + \dots + a_{n-1} z^{n-1}$$
is said to be a Littlewood polynomial if all coefficients are in $\left\{-1,1\right\}$. Balister-Bollob\'as-Morris-Sahasrabudhe-Tiba \cite{balister} showed the existence of constants $0 < c < C < \infty$ such that for each $n \in \mathbb{N}$ there exists a Littlewood polynomial satisfying
$$ c\sqrt{n} \leq |p(e^{it})| \leq C \sqrt{n}.$$
This shows that taking the coefficients of a flat Littlewood polynomial gives rise to a circulant approximate Hadamard matrix which completes the argument.
\end{proof}

\section{Remarks}

\subsection{Upper bounds.} There have been extensive efforts to find `good' polynomials: a celebrated example are the Golay-Rudin-Shapiro polynomials \cite{golay, rudin, shapiro} which are known to satisfy the upper bound $|p(e^{it})| \leq C \sqrt{n}$. Such polynomials give rise to circulants satisfying $\|Ax\| \leq C \sqrt{n} \|x\|$ and there are at least some guarantees that the lower bound is not violated on too large a subspace. We note that this is better than choosing the $\pm 1$ coefficients randomly: in that case one only gets an upper bound of $\leq C \sqrt{n \log{n}} \|x\|$. However, it is less clear whether any of this can be used to perturb the matrix $A$ into a non-circulant that is approximately Hadamard with small constants. 

\subsection{Lower bounds.} As an example in the other direction, we note the following nice and completely explicit construction by Carroll-Eustice-Figiel \cite{carroll}: if $P(z)$ is a polynomial of degree $d$ with all coefficients in $\pm 1$, then $Q(z) = P(z)P(z^{d+1})$ is a polynomial of degree $d(d+2)$ all of whose coefficients are $\pm 1$. If $|P(e^{it})| > 1$, then the minimal modulus of $Q$ is at least the square of the minimal modulus of $P$ and one obtains growth. Initializing this procedure with a very good polynomial of degree 12, Carroll-Eustice-Figiel show that the arising sequence satisfies 
$$\min |P(e^{it})| \geq (\deg P)^{0.431}.$$ This leads to a sequence of circulant matrices satisfying $\|Ax\| \geq c \cdot n^{0.431} \|x\|$ which is not as good as the main result but follows from a very simple iterative procedure. One wonders whether other such procedures might exist.

\subsection{Ultra-flat polynomials}  It is a famous open problem whether \textit{ultra-flat} polynomials exist: these are polynomials with $\pm 1$ coefficients such that 
$$(1-\varepsilon) \sqrt{n} \leq |p(e^{it})| \leq (1+\varepsilon) \sqrt{n}$$
 (for any $\varepsilon >0$ and $n$ sufficiently large depending on $\varepsilon$). We note that when one relaxes the condition $a_i = \pm 1$ and allows for complex coefficients, $|a_i| = 1$, then Kahane \cite{kah} showed that such ultraflat polynomials exist (see also Bombieri \& Bourgain \cite{bom}). This problem is also naturally connected to the problem of non-existence of Barker sequences since Barker sequences could be used to construct ultra-flat polynomials \cite{bor}. Recent extensive numerical work by Odlyzko \cite{o} suggests that ultraflat polynomials with $\pm 1$ coefficients might simply not exist. \\
 All these ideas suggest a natural conjecture which unifies a strenghtened version of Ryser's conjecture as well as a `sampling' version of the conjectured non-existence of ultra-flat polynomials. We state them separately.

\begin{conjecture}[Version A: Quantitative Ryser] There exists $\varepsilon_0 > 0$ so that for all $n > 4$ and all $n \times n$ \emph{circulants} with $\pm 1$ entries there exist a vector $x \in \mathbb{R}^n$ with
$$ \left| \|Ax\| - \sqrt{n} \cdot \|x\| \right| \geq \varepsilon_0 n^{1/4}.$$
\end{conjecture}

This conjecture says that circulant matrices with $\pm 1$ entries not only fail to be Hadamard matrices, they do so in a precise quantitative sense. We note that the power $1/4$, if true, would be optimal up to logarithmic factors: a classical argument, first given in
Matolcsi-Ruzsa \cite[Theorem 9.2]{mat} and also discussed by Jaming-Matolcsi \cite[Proposition 3.2]{jaming} and Dong-Rudelson \cite[Corollary 2.4]{dong}, shows that when $q$ is prime, a small probabilistic modification of the Legendre symbol in $\mathbb{Z}_q$ leads to circulants attaining the upper bound $n^{1/4} \sqrt{\log{n}}$ in the conjecture. The construction, while probabilistic, is explicit enough to be easily implemented on a computer (see below). It requires a `magic' ingredient, Gauss sums giving a wonderfully flat Fourier spectrum; thus, while not inconceivable that better constructions exist, one would expect that they would require `at least as much magic'.\\
By the reasoning above, one is led to a natural `sampling' variant of the ultra-flat polynomial problem: instead of asking the polynomials to be flat everywhere, one could asks them to only be flat at roots of unity. This is a very different problem.

\begin{conjecture}[Version B: Ultra-flat at Roots of Unity] There exists $\varepsilon_0 > 0$ so that for all $n > 4$ and every polynomial
$$ p(z) = a_0 + a_1 z + \dots + a_{n-1} z^{n-1} \qquad \qquad \mbox{with} \quad a_i = \pm 1$$
one has
$$ \max_{0 \leq j \leq n-1} \left|  \left| p\left(e^{\frac{2\pi i j}{n}}\right)\right| - \sqrt{n}\right| \geq \varepsilon_0 n^{1/4}.$$
\end{conjecture}
 
These two conjectures are equivalent. Since the conjecture implies Ryser's conjecture concerning the non-existence of circulant Hadamard matrices, it is presumably difficult. On the other hand, maybe it is simply too strong; that would also be interesting as it would lead to even better `almost-Hadamard' circulant matrices (for Version A) and a rather interesting sequence of polynomials (for Version B).
We note that it is already interesting that the scaling changes quite dramatically when switching from ultra-flat (from $n^{1/2}$) to ultra-flat at roots of unity (to $n^{1/4}$, at least when $n$ is prime). To illustrate this, we took the 500th prime number $n = 3571$ (with $\sqrt{n} = 59.7\dots$) and illustrated the construction in a particular (random) instance. This leads to a complex polynomial of degree $n-1$ whose behavior is shown in Fig. 1 (for a very small range of values, $p(e^{it})$ for $1 \leq t \leq 1.05$, otherwise the picture would be mainly black). 

 \begin{center}
\begin{figure}[h!]
\begin{tikzpicture}
\node at (0,0) {\includegraphics[width=0.3\textwidth]{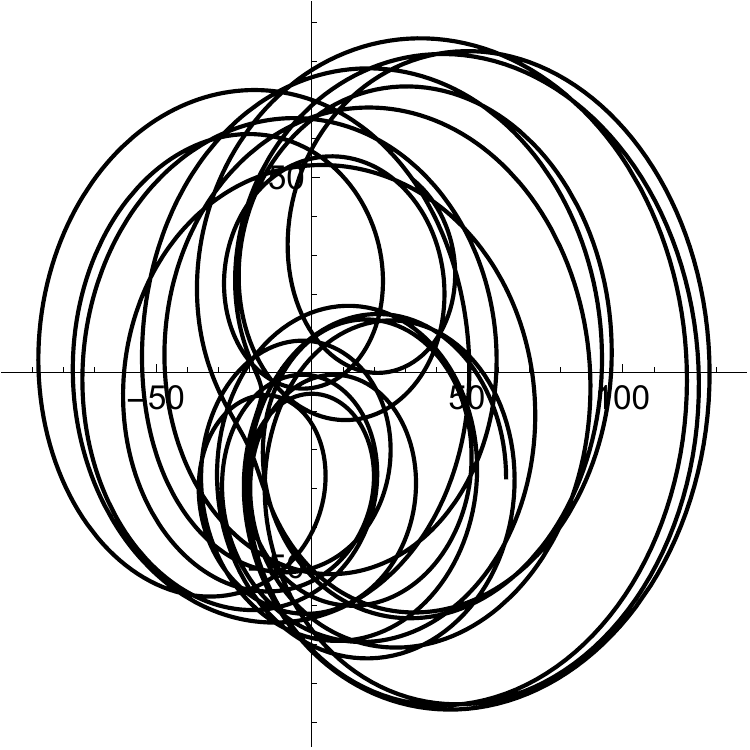}};
\node at (5,0) {\includegraphics[width=0.4\textwidth]{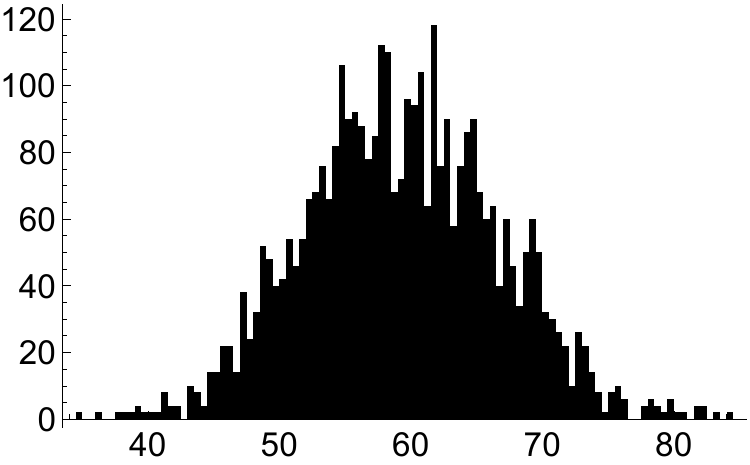}};
\end{tikzpicture}
\caption{Left: $p(e^{it})$ for $1 \leq t \leq 1.05$. Right: a histogram of the values of $| p\left(\cdot\right)|$ when evaluated at the roots of unity. }
\end{figure}
\end{center}

We see that said polynomial is far from ultra-flat, however, it has very good behavior 
when evaluated at the $n$ roots of unity. In particular, one can see hints of a Gaussian centered at $\sqrt{n}$ with standard deviation $\sim n^{1/4}$. This particular polynomial corresponds to a $3571 \times 3571$ circulant matrix with $\pm 1$ entries and with condition number $\sim 2.42$. Asymptotically, the construction leads to $n \times n$ circulant matrices with $\pm 1$ entries and a condition number of $1 + c \cdot n^{-1/4} \sqrt{\log{n}}$ ($n$ prime).


\begin{thebibliography}{10}

\bibitem{ag} S. S. Agaian, Hadamard matrices and their applications. Lecture Notes in Mathematics,
1168. Springer-Verlag, Berlin, 1985.

\bibitem{balister} P. Balister, B. Bollob\'as, R. Morris, J. Sahasrabudhe and M. Tiba, Flat Littlewood polynomials exist. Ann. Math, 192 (2020), 977-1004.

\bibitem{ban}T. Banica, I. Nechita and K. Zyczkowski, Almost Hadamard matrices: general theory and examples. Open Systems \& Information Dynamics, 19 (2012), 1250024.

\bibitem{bom}  E. Bombieri and J. Bourgain, On Kahane’s ultraflat polynomials, J. European Math. Soc., vol. 11, 2009, pp. 627--703.

\bibitem{bor} P. Borwein and M. Mossinghoff, Barker sequences and flat polynomials. London Mathematical Society Lecture note series, 352 (2008), p. 71-88.

\bibitem{carroll} F. W. Carroll, D. Eustice and T. Figiel, The minimum modulus of polynomials with coefficients of modulus one. Journal of the London Mathematical Society 2 (1977), p. 76-82.


\bibitem{golay} M. J. E. Golay, Multi-slit spectroscopy, J. Opt. Soc. Amer. 39 (1949), 437--444.


\bibitem{dong} X. Dong and M. Rudelson, Approximately Hadamard Matrices and Riesz Bases in Random Frames,
International Mathematics Research Notices (2024), p. 2044--2065


\bibitem{hor} K. J. Horadam, Hadamard Matrices and Their Applications,  Princeton University Press, 2007

\bibitem{jaming} P. Jaming and M. Matolcsi, On the existence of flat orthogonal matrices. Acta Mathematica Hungarica, 147 (2015), 179-188.

\bibitem{kah} J.-P. Kahane, Sur les polynomes a coefficients unimodulaires, Bull. London Math. Soc., vol. 12, 1980, pp. 321--342.

\bibitem{karner} H. Karner, J.  Schneid and C. Ueberhuber, Spectral decomposition of real circulant matrices. Linear Algebra and Its Applications 367 (2003), p. 301-311.

\bibitem{little}  J.E. Littlewood, On polynomials $\sum^n \pm z^m, \sum^n e^{\alpha_m i} z^m, z = e^{\theta i}$, J. London Math. Soc., 41 (1966), p.367--376.

\bibitem{mat} M. Matolcsi and I. Z. Ruzsa, Difference sets and positive exponential sums I.
General properties, Journal of Fourier Anal. Appl., 20, 17--41, (2014).

\bibitem{o} A. Odlyzko, Search for ultraflat polynomials with plus and minus one coefficients, in: Connections in Discrete Mathematics: A Celebration of the Work of Ron Graham, Cambridge University Press, 2018.


\bibitem{paley} R. Paley, On orthogonal matrices, Jour. of Mathematics and Physics 12 (1933): p. 311--320.

\bibitem{park} K. H. Park and H. Y. Song, Quasi-hadamard matrix. In 2010 IEEE International Symposium on Information Theory (pp. 1243-1247). IEEE.

\bibitem{rudin} W. Rudin, Some theorems on Fourier coefficients, Proc. Amer. Math. Soc.10 (1959),  855--859.

\bibitem{ryser} H. J. Ryser: Combinatorial Mathematics. Wiley, New York (1963).

\bibitem{schm} B. Schmidt, Cyclotomic integers and finite geometry. Journal of the American Mathematical Society, 12(4), 929--952.

\bibitem{shapiro} H.S. Shapiro, Extremal problems for polynomials, Thesis for S.M. Degree, 1952, 102 pp

\bibitem{sylv} J. Sylvester, LX. Thoughts on inverse orthogonal matrices, simultaneous sign successions, and tessellated pavements in two or more colours, with applications to Newton’s
rule, ornamental tile-work, and the theory of numbers., The London, Edinburgh, and
Dublin Philosophical Magazine and Journal of Science 34 (1867), no.232 461--475 pp.


\end{thebibliography}
\end{document}